\documentclass[a4paper, 12pt, draft, reqno]{amsart}
\usepackage[final]{graphicx}
\usepackage{amsmath}
\usepackage{amssymb}
\usepackage{latexsym}
\usepackage[cmtip,matrix,arrow]{xy}
\UseComputerModernTips
\CompileMatrices

\newcommand{\ZZ}{{\mathbb Z}}
\newcommand{\RR}{{\mathbb R}}
\newcommand{\NN}{{\mathbb N}}
\newcommand{\CC}{{\mathbb C}}

\DeclareMathOperator{\Hom}{Hom}

\DeclareMathOperator{\im}{im}

\DeclareMathOperator{\mn}{min}
\DeclareMathOperator{\res}{res}

\newcommand{\too}{\longrightarrow}

\begin{document}
\theoremstyle{plain}
\newtheorem{thm}{Theorem}[section]
\newtheorem{prop}[thm]{Proposition}
\newtheorem{lem}[thm]{Lemma}
\newtheorem{cor}[thm]{Corollary}
\newtheorem{conj}[thm]{Conjecture}
\newtheorem{claim}[thm]{Claim}
\theoremstyle{definition}
\newtheorem{rem}[thm]{Remark}
\newtheorem{ass}[thm]{Assumption}
\newtheorem{defn}[thm]{Definition}
\newtheorem{example}[thm]{Example}

\setlength{\parskip}{1ex}

\title[The blocks of the Brauer algebra]{A geometric characterisation
of\\ the blocks of the Brauer algebra}
 \author{Anton Cox}
 \email{A.G.Cox@city.ac.uk, M.Devisscher@city.ac.uk, P.P.Martin@city.ac.uk}
 \author{Maud De Visscher} 
 \author{Paul Martin} 
 \address{Centre for Mathematical Science\\
 City University\\
 Northampton Square\\ 
 London\\ 
 EC1V 0HB\\
 England.} 
\subjclass[2000]{Primary 20G05}
\begin{abstract}We give a geometric description of the
 blocks of the Brauer algebra $B_n(\delta)$ in characteristic zero as
orbits of the Weyl group of type $D_n$. We show how the corresponding
affine Weyl group controls the representation theory of the
Brauer algebra in positive characteristic, with orbits corresponding
to unions of blocks.
\end{abstract}
 
\maketitle
\medskip

\section{Introduction}

Classical Schur-Weyl duality relates the representation theory of the
symmetric and general linear groups by realising each as the
centraliser algebra of the action of the other on a certain tensor
space. The Brauer algebra $B_n(\delta)$ was introduced to provide a
corresponding duality for the symplectic and orthogonal groups
\cite{brauer}. The abstract $k$-algebra is defined for each $\delta\in
k$, however for Brauer the key case is $k=\CC$ with $\delta$ integral,
when the action of $B_n(\delta)$ on $(\CC^{|\delta|})^{\otimes n}$ can
be identified with the centraliser algebra for the corresponding group
action of O$(\delta,\CC)$ for $\delta$ positive, and with
Sp$(-\delta,\CC)$ for $\delta$ negative). In characteristic $p$, the
natural algebra in correspondence to the centraliser algebra for
$\delta$ negative is the symplectic Schur algebra
\cite{dongf,dotpoly1,oe1,ddh}.

For $|\delta|<n$ the centraliser algebra is a proper quotient of the
Brauer algebra. Thus, despite the fact that the symplectic and
orthogonal groups, and hence the centraliser, are semisimple over
$\CC$, the Brauer algebra can have a non-trivial
cohomological structure in such cases.

Brown \cite{brownbrauer} showed that the Brauer algebra is semisimple
over $\CC$ for generic values of $\delta$. Wenzl proved that
$B_n(\delta)$ is semisimple over $\CC$ for all non-integer $\delta$
\cite{wenzlbrauer}. It was not until very recently that any progress
was made in positive characteristic. A necessary and sufficient
condition for semisimplicity (valid over an arbitrary field) was given
by Rui \cite{ruibrauer}. The blocks were determined in characteristic
zero \cite{cdm} by the authors.

The block result uses the theory of towers of recollement developed in
\cite{cmpx}, and built on work by Doran, Hanlon and Wales
\cite{dhw}. The approach was combinatorial, using the language of
partitions and tableaux, and depended also on a careful analysis of
the action of the symmetric group $\Sigma_n$, realised as a subalgebra
of the Brauer algebra. However, we speculated in \cite{cdm} that there
could be an alcove geometric version, in the language of algebraic Lie
theory \cite{jannewed} (despite the absence of an obvious
Lie-theoretic context) . This should replace the combinatorics of
partitions by the action of a suitable reflection group on a weight
space, so that the blocks correspond to orbits under this action. In
this paper we will give such a geometric description of this block
result.

A priori there is no specific evidence from algebraic Lie theory to
suggest that such a reflection group action will exist (beyond certain
similarities with the partition algebra case, where there is a
reflection group of infinite type $A$ \cite{mwdef}).  As already
noted, the obvious link to Lie theory (via the duality with symplectic
and orthogonal groups) in characteristic zero only corresponds to a
semisimple quotient.

Remarkably however, we will show that there is a Weyl group $W$ of type
$D$ which does control the representation theory.  To obtain a natural
action of this group, we will find that it is easier to work with the
transpose of the usual partition notation. (This is reminiscent of the
relation under Ringel duality between the combinatorics of the
symmetric and general linear groups, although we do not have a
candidate for a corresponding dual object in this case.)

Our proof of the geometric block result in characteristic $0$ is
entirely combinatorial, as we show that the action of $W$ corresponds
to the combinatorial description of blocks in \cite{cdm}. However,
having done this, it is natural to consider extending these results to
arbitrary fields. 

As the algebras and (cell) modules under consideration can all be
defined \lq integrally' (over $\ZZ[\delta]$), one might hope that some
aspects of the characteristic $0$ theory could be translated to other
characteristics by a reduction mod $p$ argument. If this were the
case then, for consistency between different values of $\delta$ which
are congruent modulo $p$, we might expect that the role of the Weyl
group would be replaced by the corresponding affine Weyl group, so
that blocks again lie within orbits.

We will extend certain basic results in \cite{dhw} to arbitrary
characteristic, and then show that orbits of the affine Weyl group do
indeed correspond to (possibly non-trivial) unions of blocks of the
Brauer algebra.

In Section \ref{Braueris} we review some basic properties of the
Brauer algebra, following \cite{cdm}.  Sections \ref{Wis} and
\ref{Waffis} review the Weyl and affine Weyl groups of type $D$, and
give a combinatorial description of their orbits on a weight space.
Using this description we prove in Section \ref{blockzero} that we can
restate the block result from \cite{cdm} using Weyl group
orbits. Section \ref{blockp} generalises certain representation
theoretic results from \cite{dhw} and \cite{cdm} to positive
characteristic, which are then used to give a necessary condition for
two weights to lie in the same block in terms of the affine Weyl
group.

In Section \ref{absect} we describe how abacus notation \cite{jk} can
be applied to the Brauer algebra, and use this to show that the orbits
of the affine Weyl group do not give a sufficient condition for two
weights to lie in the same block.

\section{The Brauer algebra}\label{Braueris}

We begin with a very brief review of the basic theory of Brauer
algebras; details can be found in \cite{cdm}. Fix an algebraically
closed field $k$ of characteristic $p\geq 0$, and some $\delta\in
k$. For $n\in\NN$ the Brauer algebra $B_n(\delta)$ can be defined in
terms of a basis of partitions of
$\{1,\ldots,n,\bar{1},\ldots,\bar{n}\}$ into pairs. To determine the
product $AB$ of two basis elements, represent each by a graph on $2n$
points, and identify the vertices $\bar{1},\ldots,\bar{n}$ of $A$ with the
vertices $1,\ldots n$ of $B$ respectively. The graph thus obtained may
contain some number ($t$ say) of closed loops; the product $AB$ is
then defined to be $\delta^tC$, where $C$ is the basis element
corresponding to the graph arising after these closed loops are
removed (ignoring intermediate vertices in connected components).

\begin{figure}[ht]
\includegraphics{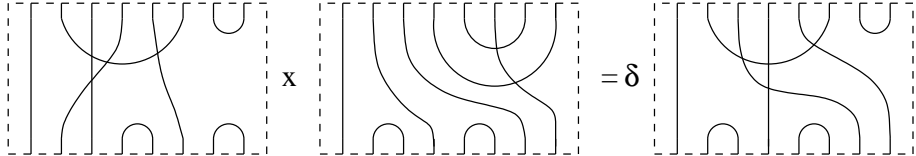}
\caption{The product of two diagrams in $B_n(\delta)$.}
\label{multex}
\end{figure}

Usually we represent basis elements graphically by a diagram with $n$
northern nodes numbered $1$ to $n$ from left to right, and $n$
southern nodes numbered $\bar{1}$ to $\bar{n}$ from left to right,
where each node is connected to precisely one other by a line. Edges
joining northern nodes to southern nodes of a diagram are called
propagating lines, the remainder are called northern or southern
arcs. An example of the product of two diagrams in given in Figure \ref{multex}.

\begin{figure}[ht]
\includegraphics{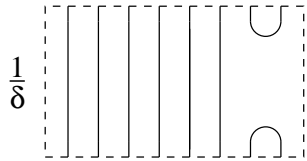}
\caption{The element $e_8$.}
\label{eis}
\end{figure}

With this convention, and assuming that $\delta\neq 0$, we have for
each $n\geq 2$ an idempotent $e_n$ as illustrated in Figure \ref{eis}.

 These idempotents induce algebra
isomorphisms
\begin{equation}\label{A1}
\Phi_n:B_{n-2}(\delta)\longrightarrow e_nB_n(\delta)e_n
\end{equation}
 which take
a diagram in $B_{n-2}$ to the diagram in $B_n$ obtained by adding an
extra northern and southern arc to the righthand end. From this we
obtain, following \cite{green}, an exact
localization functor
\begin{eqnarray*}
F_n \, :\, B_n(\delta)\mbox{\rm -mod} &\longrightarrow& 
B_{n-2}(\delta)\mbox{\rm -mod}\\
 M &\longmapsto& e_n M
 \end{eqnarray*}
and a right exact globalization functor
 \begin{eqnarray*}
 G_n\, : \, B_n(\delta)\mbox{\rm -mod} &\longrightarrow&
B_{n+2}(\delta)\mbox{\rm -mod}\\ M &\longmapsto& B_{n+2}e_{n+2}\otimes_{B_n}M.
\end{eqnarray*}
 Note that $F_{n+2}G_n(M)\cong M$ for all $M\in B_n\mbox{\rm -mod}$, and
 hence $G_n$ is a full embedding. As 
$$B_n(\delta)/B_n(\delta)e_nB_n(\delta)\cong k\Sigma_n$$ the group
algebra of the symmetric group on $n$ symbols, it follows from
\cite{green} and (\ref{A1}) that the simple $B_n$-modules are indexed
by the set
\begin{equation}\label{indexset}
\Lambda_n = \Lambda^n \sqcup \Lambda_{n-2}= \Lambda^n \sqcup
\Lambda^{n-2} \sqcup \cdots \sqcup \Lambda^{\mn}
\end{equation}
 where $\Lambda^n$ denotes an indexing set for the simple
$k\Sigma_n$-modules, and $\min=0$ or $1$ depending on the parity of
$n$. (When $\delta=0$ a slight modification of this construction is
needed; see \cite{hpbrauer} or \cite[Section 8]{cdm}.) If $\delta\neq
0$ and either $p=0$ or $p>n$ then the set $\Lambda^n$ corresponds to
the set of partitions of $n$; we write $\lambda \vdash n$ if $\lambda$
is such a partition.

If $\delta\neq 0$ and $p=0$ or $p>n$ then the algebra $B_n(\delta)$ is
quasihereditary -- in general however it is only cellular. In all
cases however we can explicitly construct a standard/cell module
$\Delta_n(\lambda)$ for each partition $\lambda$ of $m$ where $m\leq
n$ with $m-n$ even (by arguing as in \cite[Section 2]{dhw}). In the
quasihereditary case, the heads $L_n(\lambda)$ of the standard modules
$\Delta_n(\lambda)$ are simple, and provide a full set of simple
$B_n(\delta)$-modules. In the general cellular case, a proper subset
of the heads of the cell modules is sufficient to provide such a full
set of simples. The key result which we will need is that in all
cases, the blocks of the algebra correspond to the equivalence classes
of simple modules generated by the relation of occurring in the same
cell or standard module \cite[(3.9) Remarks]{gl}.

\section{Orbits of weights for the Weyl group of type $D$}
\label{Wis}

We review some basic results about the Weyl group of type $D$,
following \cite[Plate IV]{bourbaki}.  Let $\epsilon_i$ with $1\leq i\leq n$
be a set of formal symbols. We set
$$X\ (=X_n)=\bigoplus_{i=1}^{n}\ZZ\epsilon_i$$ which will play the role
of a weight lattice. We denote an element
$$\lambda=\lambda_1\epsilon_1+\cdots+\lambda_n\epsilon_n$$ in $X$ by
any tuple of the form $(\lambda_1,\ldots,\lambda_m)$, with $m\leq n$,
where $\lambda_i=0$ for $i>m$.  The set of dominant weights is given
by
$$X^+=\{\lambda\in X :
\lambda=\lambda_1\epsilon_1+\cdots\lambda_n\epsilon_n\ \mbox{\rm
with}\ \lambda_1\geq\cdots\geq\lambda_n\geq 0\}.$$ 
Define an inner
product on $E=X\otimes_\ZZ\RR$ by setting
$$(\epsilon_i,\epsilon_j)=\delta_{ij}$$ 
and extending by linearity.

Consider the root system of type $D_n$:
$$\Phi=\{\pm(\epsilon_i-\epsilon_j), \pm(\epsilon_i+\epsilon_j):1\leq
i<j\leq n\}$$ 
For each root $\beta\in\Phi$ we define a corresponding reflection
$s_{\beta}$ on $E$ by
\begin{equation}\label{reflect}
s_{\beta}(\lambda)=\lambda-(\lambda,\beta)\beta
\end{equation}
for all $\lambda\in E$, and let $W$ be the group generated by these
reflections.  Fix $\delta$ and define $\rho$ $(=\rho(\delta))\in E$ by
$$\rho=(-\frac{\delta}{2},-\frac{\delta}{2}-1,-\frac{\delta}{2}-2,\ldots,
-\frac{\delta}{2}-(n-1)).$$
We consider the dot action of $W$ on $E$ given by
$$w.\lambda=w(\lambda+\rho)-\rho$$ for all $w\in W$ and $\lambda\in
E$. (Note that this preserves the lattice $X$.) This is the action
which we will consider henceforth.

It will be convenient to have an explicit description of the dot
action of $W$ on $X$. Let $\Sigma_n$ denote the group of permutations
of ${\bf n}=\{1,\ldots, n\}$.  Given
$\lambda=(\lambda_1,\lambda_2,\ldots,\lambda_n)$ and
$\mu=(\mu_1,\mu_2,\ldots, \mu_n)$ in $X$, we have $\mu=w.\lambda$ for
some $w\in W$ if and only if
$$\mu_i+\rho_i=\sigma(i)(\lambda_{\pi(i)}+\rho_{\pi(i)})$$ for all
$1\leq i\leq n$ and some $\pi\in\Sigma_n$ and $\sigma:{\bf n}\too\{\pm
1\}$ with
$$d(\sigma)=|\{i:\sigma(i)=-1\}|$$
even. (See  \cite[IV.4.8]{bourbaki}.)
Thus $\mu=w.\lambda$ if and only if there exists $\pi\in\Sigma_n$ such
that for all $1\leq i\leq n$ we have either
\begin{equation}\label{firstcase}
\mu_i-i=\lambda_{\pi(i)}-\pi(i)\end{equation}
or
\begin{equation}\label{secondcase}
\mu_i+\lambda_{\pi(i)}-i-\pi(i)=\delta-2
\end{equation}
and (\ref{secondcase}) occurs only for an even number of $i$.

For example, if $n=5$ and $\lambda=(6,4,-2,3,5)$ then
$\mu=(-4,\delta,5,\delta-3,4)$ is in the same orbit under the dot
action of $W$, taking $\pi(1)=3$, $\pi(2)=5$, $\pi(3)=2$, $\pi(4)=1$,
$\pi(5)=4$, and $\sigma(i)=1$ for $i$ odd and
$\sigma(i)=-1$ for $i$ even.

\begin{figure}[ht]
\includegraphics{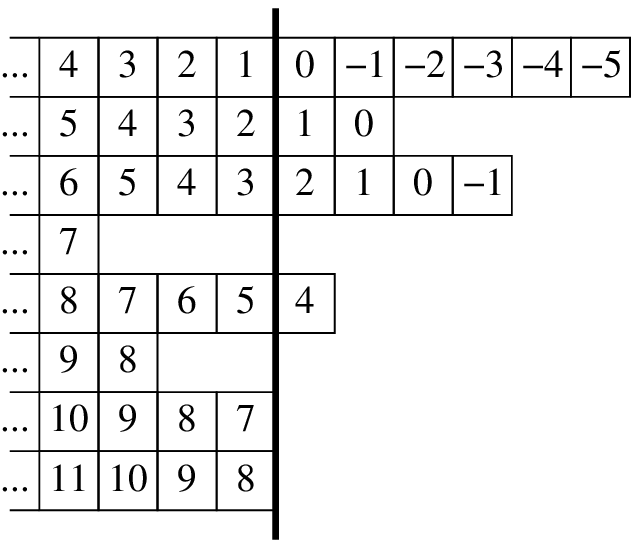}
\caption{The element $(6,2,4,-3,1,-2)$ when $n=8$.}
\label{howrep}
\end{figure}

We will also need to have a graphical representation of elements of
$X$, generalising the usual partition notation. We will represent any
$\lambda=(\lambda_1,\ldots,\lambda_n)\in X$ by a sequence of $n$ rows of
boxes, where row $i$ contains all boxes to the left of column
$\lambda_i$ inclusive, together with a vertical bar between columns
$0$ and $1$. We set the content of a box $\epsilon$ in row $i$ and
column $j$ to be $c(\epsilon)=i-j$. (This is not the usual choice for
partitions, for reasons which will become apparent later.)  For
example, when $n=8$ the element $(6,2,4,-3,1,-2)$ (and the content of
its boxes) is illustrated in Figure \ref{howrep}.

When $\lambda$ is a partition we will usually omit the portion of the
diagram to the left of the bar, and below the final non-zero row, thus
recovering the usual Young diagram notation for partitions.

If $\lambda=(\lambda_1,\ldots,\lambda_n)$ then the content
$c(\lambda)_i$ of the last box in row $i$ is $-\lambda_i+i$. Combining
this with (\ref{firstcase}) and (\ref{secondcase}) we obtain

\begin{prop}\label{orbits1}
For any two elements $\lambda$ and $\mu$ in $X$ there exists $w\in W$
with $\mu=w.\lambda$ if and only if there exists $\pi\in\Sigma_n$ and
$\sigma\, : \, {\bf n} \rightarrow \{\pm 1\}$ with $d(\sigma)$ even
such that for all $1\leq i\leq n$ we have either
$$\sigma(i)=1 \quad \mbox{and} \quad c(\mu)_i=c(\lambda)_{\pi(i)}$$
or 
$$\sigma(i)=-1 \quad \mbox{and} \quad c(\mu)_i+c(\lambda)_{\pi(i)}=2-\delta .$$
\end{prop}

It is helpful when considering low rank examples in Lie theory to use
a graphical representation of the action of a Weyl group.  As our
weight space is generally greater than two-dimensional, we can rarely 
use such an approach directly. However, we can still apply a limited
version of this approach, by considering various two-dimensional
projections of the weight lattice.

We can depict elements of the weight lattice $X$ by projecting into
the $ij$ plane for various choices of $i<j$. Each weight $\lambda$ is
represented by the projected coordinate pair $(\lambda_i,\lambda_j)$,
and each such pair represents a fibre of weights, which may
or may not include any dominant weights. For example, the point
$(0,0)$ in the $1j$ plane represents precisely one dominant weight
(the zero weight), while the $(0,0)$ point in the $23$ plane
represents the set of dominant weights
$(\lambda_1,0,0,\ldots,0)$. Clearly a necessary condition for dominance
is that $\lambda_i\geq \lambda_j\geq 0$.

\begin{figure}[ht]
\includegraphics{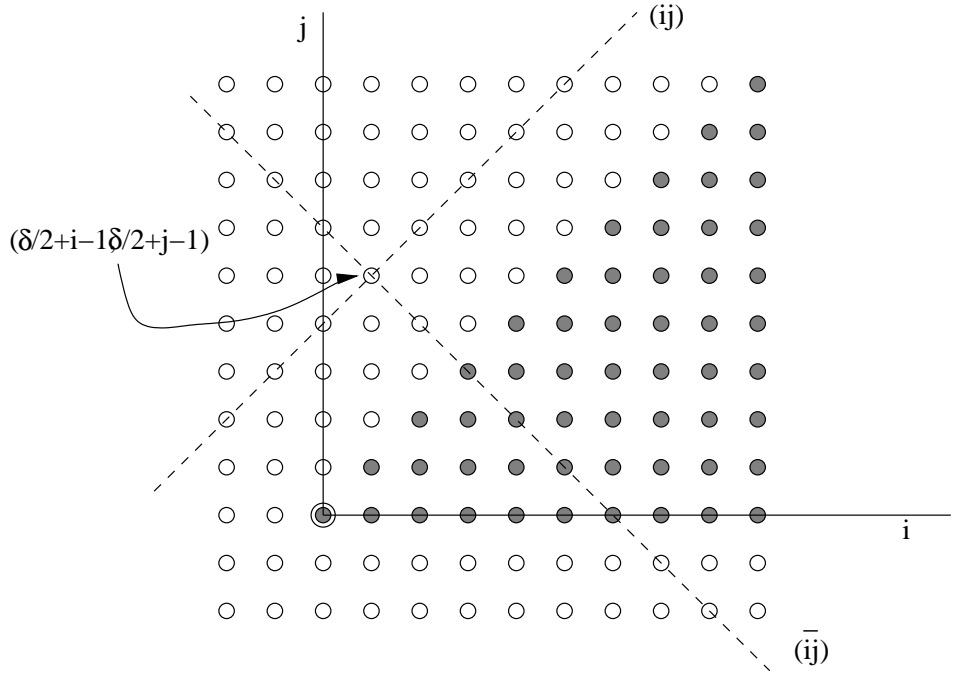}
\caption{A projection onto the $ij$ plane.}
\label{nonaffsetup}
\end{figure}

We will represent such projections in the natural two-dimensional
coordinate system, so that the set of points representing at least one
dominant weight  correspond to those shaded in Figure
\ref{nonaffsetup}. (If $\delta=2$ then the example shown is the case
$i=1$ and $j=5$.)

It will be convenient to give an explicit description of the action of
$s_{\epsilon_{i}-\epsilon_{j}}$ and $s_{\epsilon_{i}+\epsilon_{j}}$ on
a partition $\lambda$. We have that
$$s_{\epsilon_{i}-\epsilon_{j}}.\lambda=\lambda-(\lambda_i-\lambda_j-i+j)
(\epsilon_i-\epsilon_j)$$ and hence if $r=\lambda_i-\lambda_j-i+j$ is
positive (respectively negative) the effect of the dot action of
$s_{\epsilon_{i}-\epsilon_{j}}$ on $\lambda$ is to add $r$ boxes to row $j$
(respectively row $i$) and subtract $r$ boxes from row $i$
(respectively row $j$). Similarly, 
$$s_{\epsilon_{i}+\epsilon_{j}}.\lambda
=\lambda-(\lambda_i+\lambda_j-\delta+2-i-j) (\epsilon_i+\epsilon_j)$$
and hence if $r=\lambda_i+\lambda_j-\delta+2-i-j$ is positive
(respectively negative) the effect of the dot action of
$s_{\epsilon_{i}+\epsilon_{j}}$ on $\lambda$ is to remove
(respectively add) $r$ boxes from each of rows $i$ and $j$. In terms
of our projection onto the $ij$ plane these operations correspond to
reflection about the dashed lines in Figure \ref{nonaffsetup} labelled
$(ij)$ for $s_{\epsilon_{i}-\epsilon_{j}}$ and $(\bar{ij})$ for
$s_{\epsilon_{i}+\epsilon_{j}}$. Note that the position of
$(\bar{ij})$ depends on $\delta$, but $(ij)$ does not.

\begin{figure}[ht]
\includegraphics{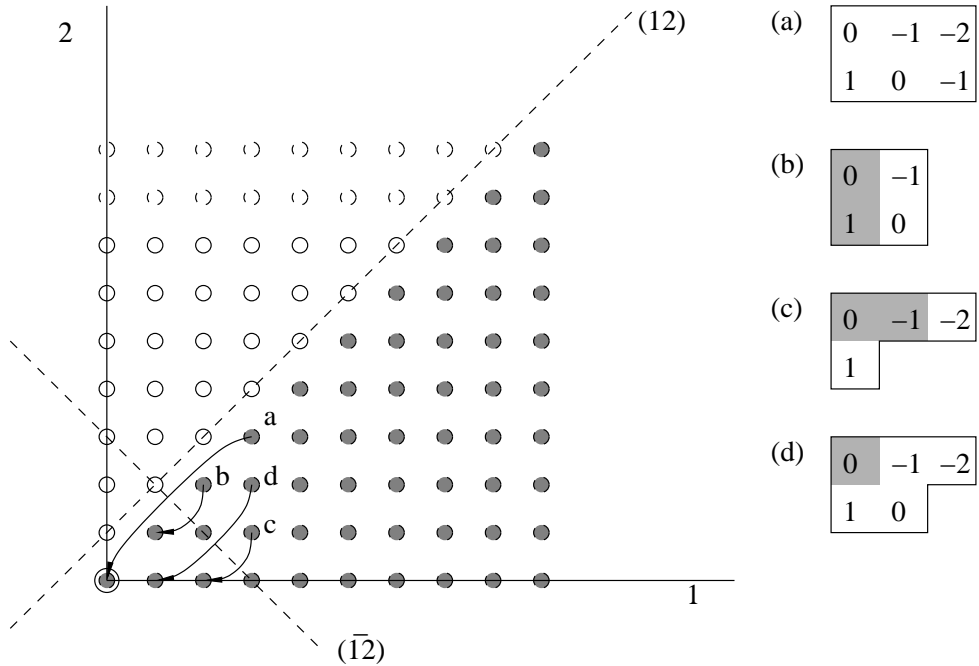}
\caption{Projections into the $12$ plane with $\delta=2$.}
\label{nona}
\end{figure}

\begin{figure}[ht]
\includegraphics{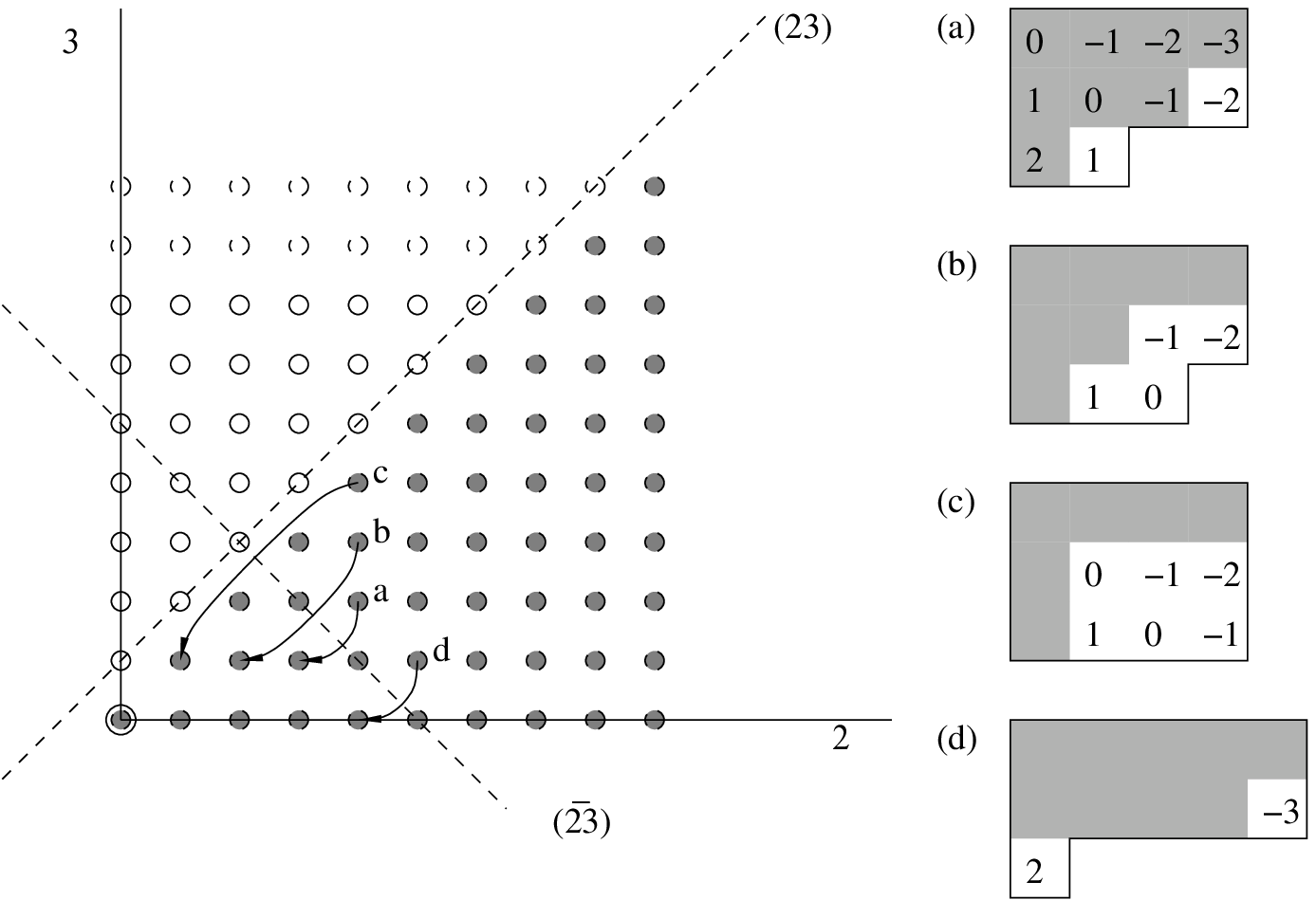}
\caption{Projections into the $23$ plane with $\delta=2$.}
\label{nonb}
\end{figure}

\begin{figure}[ht]
\includegraphics{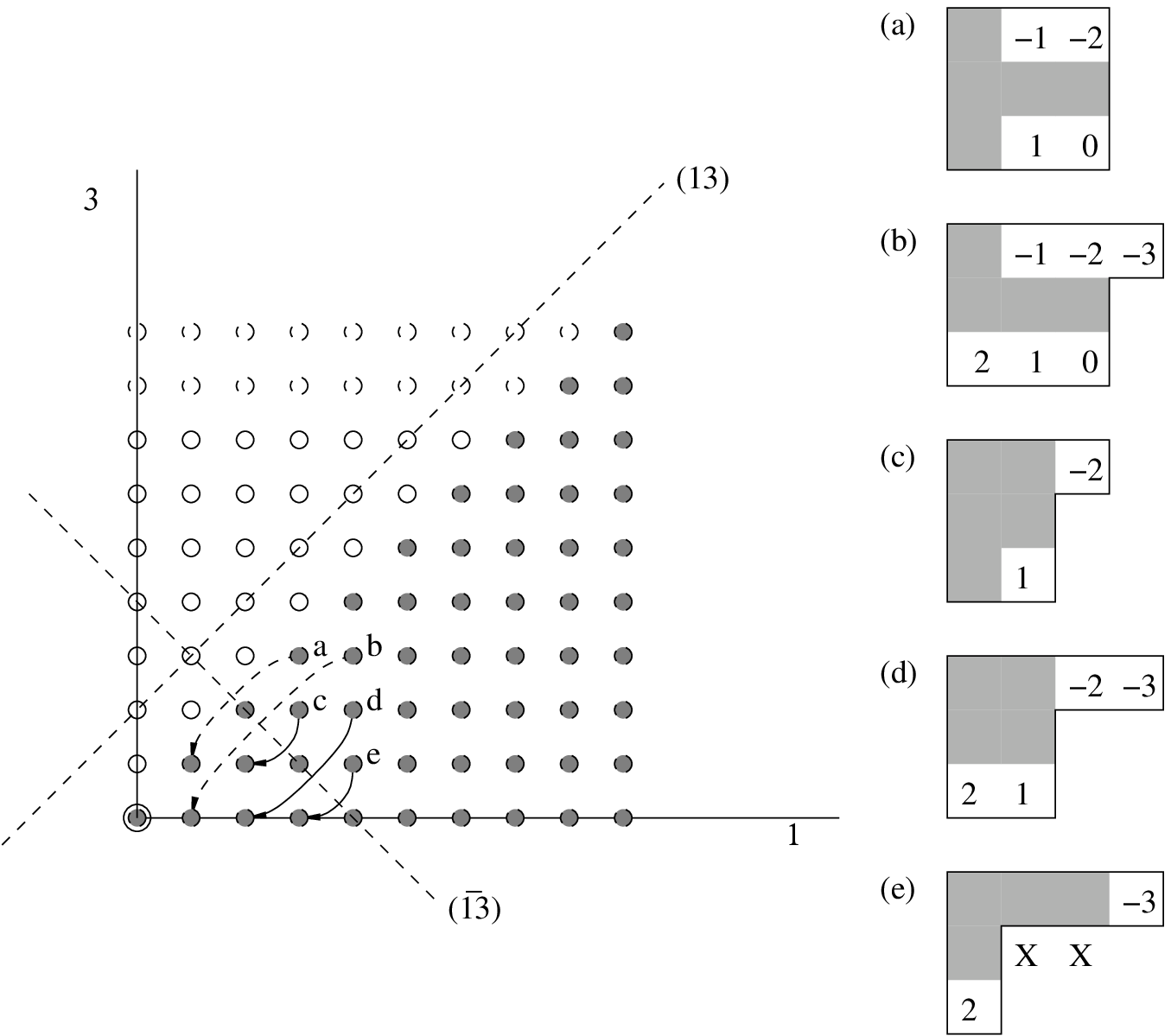}
\caption{Projections into the $13$ plane with $\delta=2$.}
\label{nonc}
\end{figure}

Various examples of reflections, and their effect on a dominant
representative of each coordinate pair, are given in Figures 
\ref{nona}, \ref{nonb}, and \ref{nonc}. For each reflection indicated,
a dominant weight is illustrated, together with a shaded
subcomposition corresponding to the image of that weight under the
reflection. Where no shading is shown (as in Figure \ref{nona}(a)) the
image is the empty partition.

Note that some of the reflections may take a dominant weight to a
non-dominant one, even if the associated fibres both contain dominant
weights. For example the cases in Figure \ref{nonc}(a) and (b)
correspond to the reflection of $(3,3,3)$ to $(1,3,1)$ and of
$(4,3,3)$ to $(1,3,0)$. Also, some reflections may represent a family
of reflections of dominant weights, as in Figure \ref{nonc}(e), where
there are three possible weights in each fibre (corresponding to
whether none, one or both of the boxes marked X are included).

\section{The blocks of the Brauer algebra in characteristic zero}
\label{blockzero}

The main result in \cite{cdm} was the determination of the blocks of
$B_n(\delta)$ when $k=\CC$. In that paper, the blocks were described
by a combinatorial condition on partitions. We would like to have a
geometric formulation of this result.

\begin{figure}[ht]
\includegraphics{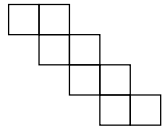}
\caption{A potentially unbalanced configuration.}
\label{exclude}
\end{figure}

We will identify the simple $B_n(\delta)$-modules with weights in
$X^+$ using the correspondence
$$(\lambda \in X^+)\quad \longleftrightarrow \quad L(\lambda^T)$$
where $\lambda^T$ denotes the conjugate partition of $\lambda$
(i.e. the one obtained by reversing the roles of rows and columns in
the usual Young diagram).  Using this correspondence, we restate the
main result of \cite{cdm} as follows.  Given two partitions
$\mu\subset\lambda$ we write $\lambda/\mu$ for the associated skew
partition. We say that that two weights $\lambda , \mu\in X^+$ are
($\delta$-)\textit{balanced} if and only if the boxes of $\lambda
/(\lambda \cap \mu)$ (respectively $\mu / (\lambda \cap \mu)$) can be
paired up such that the contents of each pair sum to $1-\delta$ and if
$\delta$ is even and the boxes with content $-\frac{\delta}{2}$ and
$\frac{2-\delta}{2}$ are configured as in Figure \ref{exclude} then
the number of rows in Figure \ref{exclude} must be even.

Noting that the definition of content given in Section 1 is the
transpose of the one used in \cite{cdm}, it is easy to see (simply by
transposing everything) that \cite{cdm} Corollary 6.7 becomes

\begin{thm}\label{oldblock}
Two simple $B_n(\delta)$-modules $L(\lambda^T)$ and $L(\mu^T)$ are in
the same block if and only if $\lambda$ and $\mu$ are balanced.
\end{thm}

We now give the desired geometric formulation of Theorem \ref{oldblock}.

\begin{thm}\label{geomblock}
Two simple $B_n(\delta)$-modules $L(\lambda^T)$ and $L(\mu^T)$ are in
the same block if and only if $\mu\in W\cdotp \lambda$
\end{thm}


\begin{proof}
We will show that this description is equivalent to that given in
Theorem \ref{oldblock}, by proceeding in two stages. First we will
show, using the action of the generators of $W$ on $X$, that two
partitions in the same orbit are balanced. This implies that the
blocks are unions of $W$-orbits. Next we will show that balanced
partitions lie in the same $W$-orbit.

\noindent
{\bf Stage 1:} The case $n=2$ is an easy calculation. For $n>2$, note
that
$$s_{\epsilon_i-\epsilon_j}
=s_{\epsilon_j+\epsilon_k}s_{\epsilon_i+\epsilon_k}
s_{\epsilon_j+\epsilon_k}$$ where $i\neq k\neq j$, and so $W$ is
generated by reflections of the form $s_{\epsilon_i +\epsilon_j}$. Now
consider the action of such a generator on a weight in $X$.
$$s_{\epsilon_i+\epsilon_j}.\lambda
=\lambda-(\lambda_i+\lambda_j-\delta+2-i-j)(\epsilon_i+\epsilon_j).$$
If $\lambda_i+\lambda_j-\delta+2-i-j\geq 0$ then this involves the
removal of two rows of boxes with respective contents
$$-\lambda_i+i+\lambda_i+\lambda_j-i-j-\delta+1,\ldots,-\lambda_i+i+1,
-\lambda_i+i$$
and
$$-\lambda_j+j+\lambda_i+\lambda_j-i-j-\delta+1,\ldots,-\lambda_j+j+1,
-\lambda_j+j$$
which simplify to 
$$\lambda_j-j-\delta+1,\ldots,-\lambda_i+i+1,
-\lambda_i+i$$
and
$$\lambda_i-i-\delta+1,\ldots,-\lambda_j+j+1, -\lambda_j+j.$$ If we
pair these two rows in reverse order, each pair of contents sum to
$1-\delta$. Note also that for $\delta$ even, the number of horizontal
pairs of boxes of content $-\frac{\delta}{2}$ and $\frac{2-\delta}{2}$
is either unchanged or decreased by 2.  The argument when
$\lambda_i-\lambda_j-\delta+2-i-j< 0$ is similar (here we add paired
boxes instead of removing them).

Now take two partitions $\lambda , \mu\in X^+$ with $\mu=w\cdotp
\lambda$ for some $w\in W$. We need to show that they are balanced,
i.e. that the boxes of $\lambda / \lambda \cap \mu$ (respectively $\mu
/ \lambda \cap \mu$) can be paired up in the appropriate way. First
observe that the set of contents of boxes in $\lambda / \lambda \cap
\mu$ and in $\mu / \lambda \cap \mu$ are disjoint. To see this,
suppose that there is a box $\epsilon$ in $\lambda / \lambda \cap \mu$
with the same content as a box $\eta$ in $\mu / \lambda \cap
\mu$. Then these two boxes must lie on the same diagonal. Assume,
without loss of generality, that $\epsilon$ appears in an earlier row
than $\eta$. As $\eta$ belongs to $\mu$ and $\epsilon$ is above and to
the left of $\eta$, we must have that $\epsilon$ is also in $\mu$ (as
$\mu$ is a partition). But then $\epsilon$ belongs to $\lambda \cap
\mu$ which is a contradiction.

Let us concentrate on the action of $w$ on boxes either with a fixed
content $c$ say or with the paired content $1-\delta -c$.  As $w$ can
be written as a product of the generators considered above, it will add and
remove pairs of boxes of these content, say
$$(\tau_1 + \tau_1')+(\tau_2 + \tau_2')+ \ldots +  (\tau_m + \tau_m')$$
$$-(\sigma_1 + \sigma_1')-(\sigma_2 + \sigma_2') - \ldots -(\sigma_q +
\sigma_q') $$ for some boxes $\tau_i$, $\tau_i'$, $\sigma_j$ and
$\sigma_j'$ with $c(\tau_i)=c=1-\delta - c(\tau_i')$ for $1\leq i \leq
m$ and $c(\sigma_j)=c=1-\delta -c(\sigma_j')$ for $1\leq j \leq q$.
Thus the number of boxes in $\mu=w\cdotp \lambda$ of content $c$
(resp. $1-\delta -c$) minus the number of boxes in $\lambda$ of
content $c$ (resp. $1-\delta -c$) is equal to $m-q$. But this must be
equal to the number of boxes in $\mu / (\lambda \cap \mu)$ of content
$c$ (resp. $1-\delta -c$) minus the number of boxes in $\lambda /
(\lambda \cap \mu)$ of content $c$ (resp. $1-\delta -c$).  As we have
just observed that the contents of boxes in $\lambda / (\lambda \cap
\mu)$ and in $\mu / (\lambda \cap \mu)$ are disjoint, we either have
$m-q \geq 0$ and
\begin{eqnarray*}
m-q &=& |\{ \mbox{boxes of content $c$ in $\mu / (\lambda \cap
\mu)$}\} |\\ &=& | \{ \mbox{boxes of content $1-\delta -c$ in $\mu /
(\lambda \cap \mu)$}\} |
\end{eqnarray*}
or $m-q<0$ and 
\begin{eqnarray*}
m-q &=& - | \{ \mbox{boxes of content $c$ in $\lambda / (\lambda \cap
\mu)$}\} |\\ &=& - | \{ \mbox{boxes of content $1-\delta -c$ in
$\lambda / (\lambda \cap \mu)$}\} |
\end{eqnarray*}
Thus the boxes of $\lambda / (\lambda \cap \mu)$ (resp. $\mu /
(\lambda \cap \mu)$) can be paired up such that the sum of the
contents in each pair is equal to $1-\delta$. Moreover, for $\delta$
even, as each generator $s_{\epsilon_i + \epsilon_j}$ either adds or
removes 2 (or no) horizontal pairs of boxes of contents
$-\frac{\delta}{2}$ and $\frac{2-\delta}{2}$, we see that $\lambda$
and $\mu$ are indeed balanced.

\noindent{\bf Stage 2:} We need to show that if $\lambda$ and $\mu$
are balanced partitions then they are in the same $W$-orbit. Note that
if $\lambda$ and $\mu$ are balanced then by definition so are
$\lambda$ and $\lambda \cap \mu$, and $\mu$ and $\lambda \cap
\mu$. Thus it is enough to show that if $\mu \subset \lambda$ are
balanced then they are in the same $W$-orbit.

We will show that whenever we have a weight $\eta \in X$ with $\eta +
\rho \in X^+$ such that $\mu \subset \eta$ (i.e. $\mu_i \leq \eta_i$ for all
$i$) are balanced, we can construct $\eta^{(1)}\in W_p \cdotp \eta$
such that either $\eta^1=\mu$ or $\eta ^{(1)}\subset \eta$ having the
same properties as $\eta$.  Starting with $\eta = \lambda$ and applying
induction will prove that $\mu \in W_p\cdotp \lambda$.

Pick a box $\epsilon$ in $\eta / \mu$ such that

(i) it is the last box in a row of $\eta$,

(ii) $\frac{1-\delta}{2} -c(\epsilon)$ is maximal.

\noindent If more than one such box exists, pick the southeastern-most
one. Say that $\epsilon$ is in row $i$. Find a box $\epsilon '$ on the
edge of $\eta / \mu$ (i.e. a box in $\eta/\mu$ such that there is no
box to the northeast, east, or southeast of it in $\eta/\mu$) with
$c(\epsilon)+c(\epsilon ')=1-\delta$. Say that $\epsilon '$ is in row
$j$.

Note that $i\neq j$ as if $i$ were equal to $j$ then there would
either be a box of content $\frac{1-\delta}{2}$ (for $\delta$ odd) or
a pair of boxes of content $-\frac{\delta}{2}$ and
$\frac{2-\delta}{2}$ (for $\delta$ even) in between $\epsilon$ and
$\epsilon '$. Now, as $\eta / \mu$ is balanced and
$\eta_{i-1}-\eta_i\geq -1$, it must contain another such box or pair of
boxes of the same content(s) in row $i-1$, as illustrated in Figure
\ref{2row1} (where the shaded area is part of $\mu$). As $\eta
_{i-1}-\eta_i\geq -1$ we see that $\eta /\mu$ contains at least two
boxes of content $c(\epsilon)$. But as $\epsilon$ was chosen with
maximal content and $\mu$ is a partition, $\eta / \mu$ can only have
one box of content $c(\epsilon ')$, as otherwise the box $\chi$ would
be in $\eta/\mu$. This contradicts the fact that
$\eta / \mu$ is balanced.

\begin{figure}[ht]
\includegraphics{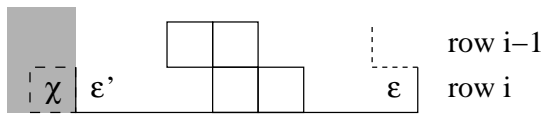}
\caption{The (impossible) configuration occurring if $i=j$.}
\label{2row1}
\end{figure}

Now let $\alpha$ be the last box in row $j$ and let $\alpha '$ be the
southeastern-most box on the edge of $\eta / \mu$ having content
$c(\alpha)=1-\delta - c(\alpha)$. Say that $\alpha '$ is in row $k$.

\noindent \textit{Case 1:} $k=j$.\\ In this case there must either be
a box of content $\frac{1-\delta}{2}$ (for $\delta$ odd) or a pair of
boxes of content $-\frac{\delta}{2}$ and $\frac{2-\delta}{2}$ (for
$\delta$ even) in between $\alpha '$ and $\alpha$. Now, as $\eta /
\mu$ is balanced and $\eta_{j-1}-\eta_j\geq -1$, it must contain
another such box or pair of boxes of the same content(s) in row $j-1$,
as illustrated in Figure \ref{2row2}.

\begin{figure}[ht]
\includegraphics{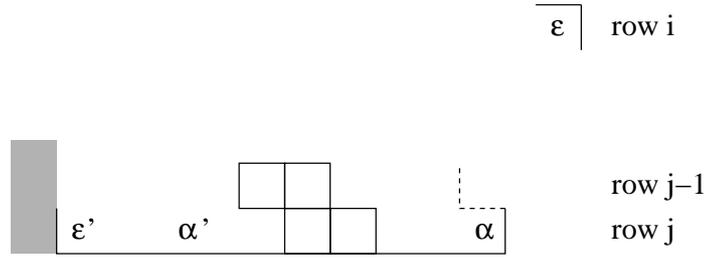}
\caption{The case $j=k$.}
\label{2row2}
\end{figure}

For each $c(\epsilon )\leq c <c(\alpha)$, define $i_c$ by saying that
the southeastern-most box of content $c$ on the edge of $\eta / \mu$
is in row $i_c$. For $c=c(\alpha)$, define $i_{c(\alpha)}=j-1$. Note
that the $i_c$'s are not necessarily all distinct. Consider all
distinct values of $i_c$ and order them
$$i=i_{c_0}<i_{c_1}<\ldots < i_{c_l}=j-1.$$
Now consider
$$\eta^{(1)} = (s_{\epsilon_i -\epsilon_{i_{c_1}}}\ldots
s_{\epsilon_i-\epsilon_{i_{c_{l-1}}}} s_{\epsilon_i -\epsilon_{j-1}}
s_{\epsilon_i+\epsilon_j})\cdotp \eta.$$ This is illustrated
schematically in Figure \ref{2row3}, where curved lines indicate
boundaries whose precise configuration does not concern us. Then $\mu
\subset \eta^{(1)}\subset \eta$ are balanced and $\eta^{(1)}+\rho \in
X^+$ as required.

\begin{figure}[ht]
\includegraphics{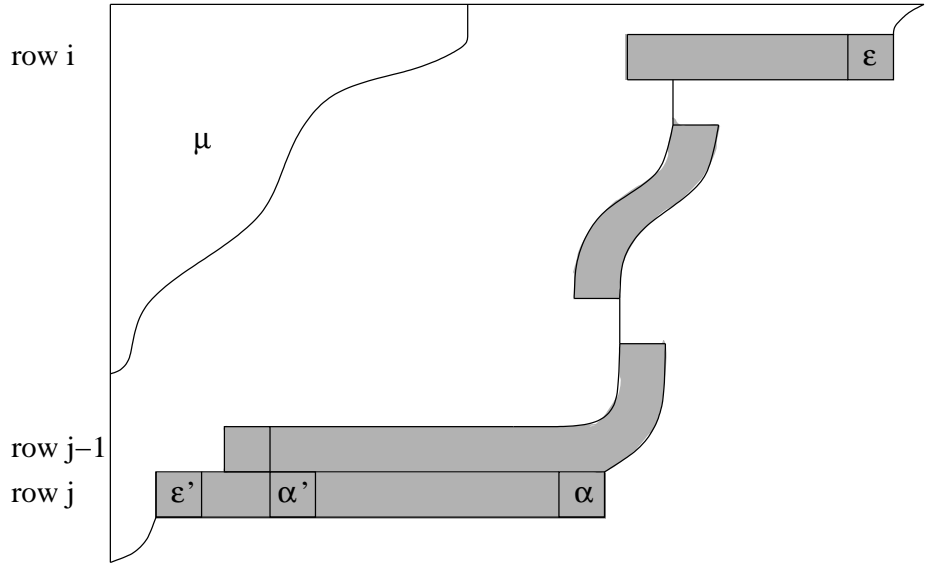}
\caption{The elements $\mu\subset\eta^{(1)}\subset\eta$.}
\label{2row3}
\end{figure}

\noindent \textit{Case 2:} $k\neq j$.\\
If $i=k$ then consider 
$$\eta^{(1)}=s_{\epsilon_i+\epsilon_j}\cdotp \eta$$ then $\mu \subset
\eta^{(1)} \subset \eta$ are balanced and $\eta^{(1)}+\rho \in X^+$.

\begin{figure}[ht]
\includegraphics{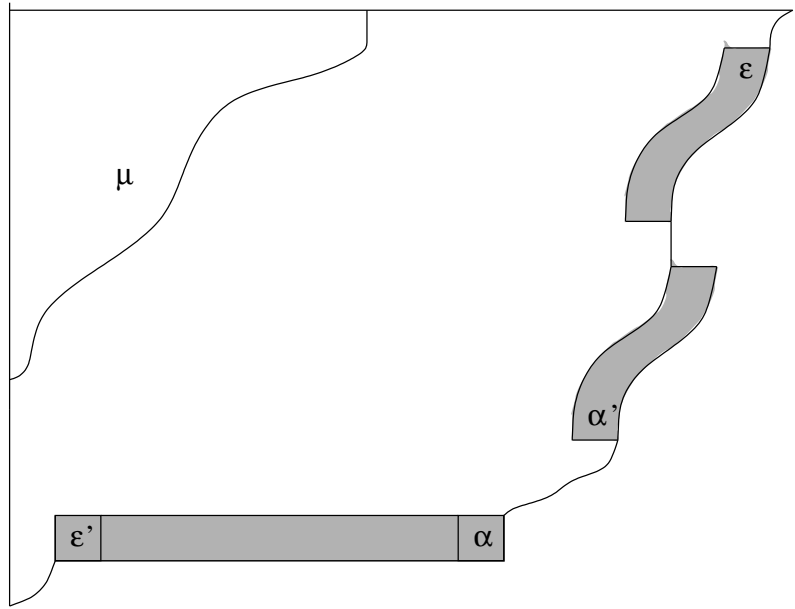}
\caption{The elements $\mu\subset\eta^{(1)}\subset\eta$.}
\label{2row4}
\end{figure}

If $k\neq i$, then as in Case 1, for each $c(\epsilon)\leq c\leq
c(\alpha ')$ we define $i_c$ by saying that the southeastern most box
in $\eta/\mu$ is in row $i_c$. As before, there are not necessarily
all distinct but we can pick a set of representatives
$$i=i_{c_0} < i_{c_1} < \ldots < i_{c_l}=k.$$
Now consider 
$$\eta^{(1)} = (s_{\epsilon_i -\epsilon_{i_{c_1}}}\ldots
s_{\epsilon_i-\epsilon_{i_{c_{l-1}}}} s_{\epsilon_i -\epsilon_{k}}
s_{\epsilon_i+\epsilon_j})\cdotp \eta.$$ Again, this is illustrated
schematically in Figure \ref{2row4}, where curved lines indicate
boundaries whose precise configuration does not concern us.
As before we have $\mu \subset \eta^{(1)} \subset \eta$ are balanced
and $\eta^{(1)} +\rho \in X^+$.
 \end{proof} 

\begin{example} We illustrate Stage 2 of the proof by an example. 
Take $\lambda = (8,8,8,7,3,3,2)$ and $\mu = (6,5,1,1)$ and $\delta
=2$. Then it is easy to see that $\mu \subset \lambda$ are
balanced. We will construct $w\in W_p$ such that $\mu = w\cdotp
\lambda$.

\begin{figure}[ht]
\includegraphics{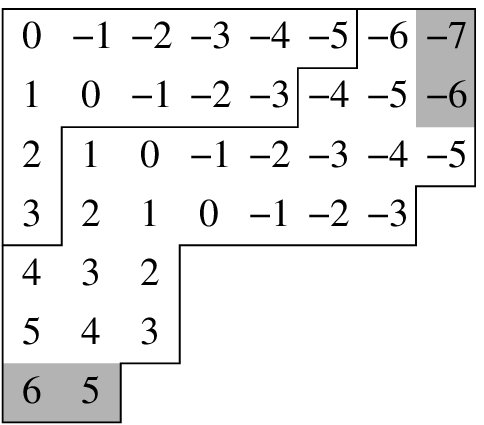}
\caption{The elements $\mu\subset\lambda^{(1)}\subset\lambda$.}
\label{ex1step1}
\end{figure} 

First consider $\lambda ^{(1)} = s_{\epsilon_1
  -\epsilon_2}s_{\epsilon_1 + \epsilon_7} \cdotp \lambda$. The
elements $\lambda$ and $\mu$ are illustrated in outline in Figure
\ref{ex1step1}, with the boxes removed to form $\lambda^{(1)}$ shaded.

Repeating the process we next consider $\lambda^{(2)}= s_{\epsilon_1 -
  \epsilon_3}s_{\epsilon_1 + \epsilon_6}\cdotp \lambda^{(1)}$, as in
  Figure \ref{ex1step2}, followed by $\lambda^{(3)}=s_{\epsilon_2 -
  \epsilon_4}s_{\epsilon_2 + \epsilon_5}\cdotp \lambda ^{(2)}$ as in
  Figure \ref{ex1step3}. Finally consider
  $\lambda^{(4)}=s_{\epsilon_3+\epsilon_4} \cdotp \lambda^{(3)}$ as
  shown in Figure \ref{ex1step4}.

\begin{figure}[ht]
\includegraphics{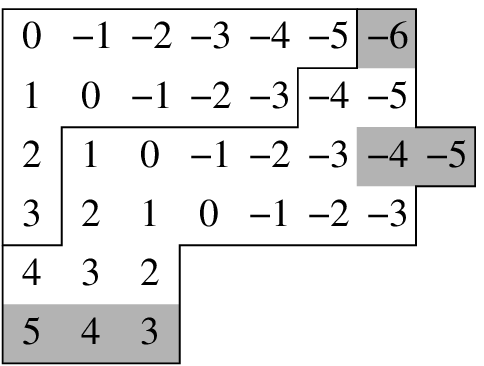}
\caption{The elements $\mu\subset\lambda^{(2)}\subset\lambda^{(1)}$.}
\label{ex1step2}
\end{figure}

\begin{figure}[ht]
\includegraphics{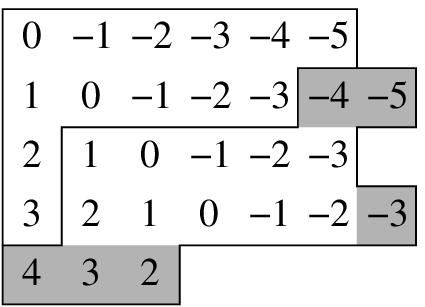}
\caption{The elements $\mu\subset\lambda^{(3)}\subset\lambda^{(2)}$.}
\label{ex1step3}
\end{figure}

\begin{figure}[ht]
\includegraphics{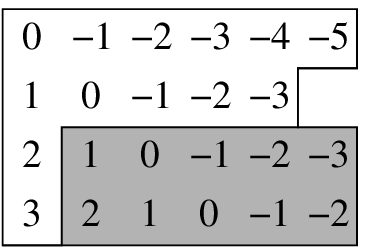}
\caption{The elements $\mu=\lambda^{(4)}\subset\lambda^{(3)}$.}
\label{ex1step4}
\end{figure} 

\end{example}

\section{Orbits of the affine Weyl group of type $D$}
\label{Waffis}

We would like to have a block result in characteristic $p>0$ similar
in spirit to Theorem \ref{geomblock}. For this we first need a
candidate to play the role of $W$. To motivate our choice of such, we
begin by considering a possible approach to modular representation
theory via reduction from characteristic $0$.
\medskip

The verification that the Brauer algebra is cellular is a
characteristic-free calculation over $\ZZ[\delta]$. Thus all of our algebras
and cell modules have a $\ZZ[\delta]$-form, from which the corresponding
objects over $k$ can be obtained by specialisation. If the maps
between cell modules that have been constructed in charateristic zero
in \cite{dhw,cdm} also had a corresponding integral form, then they
would also specialise to maps in characteristic $p$. As there is not
yet an explicit construction of these maps, we are unable to verify
this except in very small examples. However, if we assume for the
moment that it holds, this will suggest a candidate for our new
reflection group.

We will wish to consider the dot action of $W$ for different values of
shift parameter $\rho$. In such cases we will write 
$w._{\delta}\lambda$ for the element
$$w(\lambda+\rho(\delta))-\rho(\delta).$$ When we wish to emphasise
the choice of dot action we will also write $W^{\delta}$ for $W$.

Fix $\delta\in k$, and suppose that maps between cell modules in
characteristic $0$ do reduce mod $p$. Then we would expect weights to
be in the same block in characteristic $p$ if they are linked by the
action of $W^{\delta}$ in characteristic zero. However, all elements
of the form $\delta+rp$ in characteristic zero reduce to the same
element $\delta$ mod $p$, and so weights should be in the same block
if they are linked by the action of $W^{\delta+rp}$ for some
$r\in\ZZ$.  Thus our candidate for a suitable reflection group will be
${\bf W}=\langle W^{\delta+rp}: r\in\ZZ\rangle$. 

Note however that a block result does not follow \emph{automatically}
from the integrality assumption, as: (i) the chain of reflections from
${\bf W}$ linking two weights might leave the set of weights for
$B_n(\delta)$; (ii) in characteristic $p$ there may be new connections
between weights not coming from connections in characteristic zero. We
shall see that the former is indeed a problem, but that the latter
does not occur.

Now fix a prime number $p>2$ and consider the affine Weyl group
$W_p$ associated to $W$. This is defined to be
$$W_p=\langle s_{\beta,rp}:\beta\in\Phi, r\in\ZZ\rangle$$
where 
$$s_{\beta,rp}(\lambda)=\lambda
-((\lambda,\beta)-rp)\beta.$$
As before, we consider the dot action of $W_p$ on $X$ (or $E$) given by
$$w.\lambda=w(\lambda+\rho)-\rho.$$

It is an easy exercise to show

\begin{lem}\label{genrel} For all $r\in\ZZ$ and $1\leq i\neq j\neq
 k\leq n$, we have
\begin{eqnarray*} 
s_{\epsilon_i+\epsilon_j}._{\delta+rp}\lambda
&=&s_{\epsilon_i+\epsilon_j,rp}._{\delta}\lambda.\\
s_{\epsilon_i-\epsilon_j}._{\delta+rp}\lambda
&=&s_{\epsilon_i-\epsilon_j}._{\delta}\lambda.\\
s_{\epsilon_i-\epsilon_j,rp}&=& 
s_{\epsilon_j+\epsilon_k}
s_{\epsilon_i+\epsilon_k,rp}s_{\epsilon_j+\epsilon_k}.
\end{eqnarray*}
and $$s_{\epsilon_i+\epsilon_j,rp}s_{\epsilon_i+\epsilon_j}\ \mbox{
is translation by}\ rp(\epsilon_i+\epsilon_j).$$
In particular, for $n>2$ we have
$$W_p=\langle s_{\epsilon_i+\epsilon_j,rp}:
1\leq i<j\leq n\ \mbox{ and}\ r\in\ZZ\rangle.$$
\end{lem}
It follows from the first two parts of the Lemma that the group
$$W^{[r]}=\langle
s_{\epsilon_i+\epsilon_j,rp},s_{\epsilon_i-\epsilon_j} :1\leq i<j\leq
n\rangle$$ is isomorphic to the original group $W^{\delta+rp}$, and
its $\delta$-dot action on $X$ is the same as the $(\delta+rp)$-dot
action of $W^{\delta+rp}$ on $X$. Further, the usual dot action of
$W_p$ on $X$ is generated by all the $W^{[r]}$ with $r\in\ZZ$. Thus we
have

\begin{cor} For $p>2$ we have ${\bf W}\cong W_p$, and the isomorphism
  is compatible with their respective dot actions on $X$.
\end{cor}

The above considerations suggest that the affine Weyl group is a
potential candidate for the reflection group needed for a positive
characteristic block result. It will be convenient to have a
combinatorial description of the orbits of this group on $X$.

\begin{prop}\label{orbits2}
Suppose that $\lambda$ and $\mu$ are in $X$ with $|\lambda|-|\mu|$
even. Then $\mu\in W_p\cdotp \lambda$ if and only if there exists
$\pi\in\Sigma_n$ and $\sigma\, : \, {\bf n} \rightarrow \{\pm 1\}$
with $d(\sigma)$ even such that for all $1\leq i\leq
n$ we have either
$$\sigma(i)=1 \quad \mbox{and} \quad c(\mu)_i\equiv
c(\lambda)_{\pi(i)} \mod p$$ or
$$\sigma(i)=-1 \quad \mbox{and} \quad
c(\mu)_i+c(\lambda)_{\pi(i)}\equiv 2-\delta \mod p$$
\end{prop}

\begin{proof}
We have $\mu \in W_p\cdotp \lambda$ if and only if 
$$\mu+\rho=w(\lambda+\rho)+p\nu$$ for some $w\in W$ and 
$\nu\in\ZZ\Phi$. Note that for any $\nu\in X$ we have $\nu\in\ZZ\Phi$ if and
only if $|\nu|=\sum\nu_i$ is even, as 
$$2\epsilon_i=(\epsilon_i-\epsilon_{i+1})+(\epsilon_i+\epsilon_{i+1})$$
and
$$2\epsilon_{i+1}=(\epsilon_i+\epsilon_{i+1})-(\epsilon_i-\epsilon_{i+1})$$
for all $1\leq i\leq n-1$.
Thus, if $|\lambda |- |\mu|$ is even, then $\mu \in W_p\cdotp \lambda$
if and only if $\mu +\rho = w(\lambda + \rho) +p\nu$ for some $w\in W$
and some $\nu\in X$. Combining this with Proposition \ref{orbits1}
gives the result. 
\end{proof}

\begin{figure}[ht]
\includegraphics{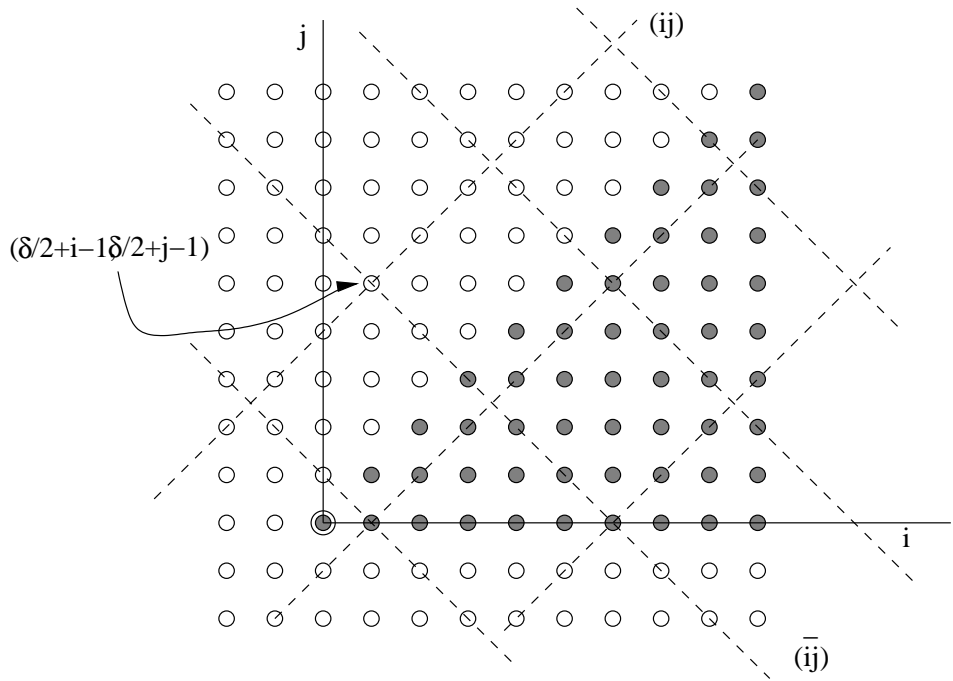}
\caption{A projection onto the $ij$ plane with $p=5$.}
\label{affsetup}
\end{figure}

As in the non-affine case, we may represent reflections graphically
via projection into the plane. In this case each projection will
contain two families of reflections; those parallel to
$s_{\epsilon_i-\epsilon_j}$ and those parallel to
$s_{\epsilon_i+\epsilon_j}$. This is illustrated for $p=5$ in Figure
\ref{affsetup}. An example of the effect of various such reflections
on partitions will be given in Figure \ref{rosetta}, after we have
introduced a third, abacus, notation.

\section{On the blocks of the Brauer algebra in characteristic $p$}
\label{blockp}

We have already seen that the blocks of the Brauer algebra in
characteristic $0$ are given by the restriction of orbits of $W$ to
the set of partitions.  We would like a corresponding result in
characteristic $p>0$ involving the orbits of $W_p$. As noted in the
introduction, one does not expect the blocks of the Brauer algebra to
be given by $W_p$ in {\it exactly} the same manner as in
characteristic $0$. Instead, we can ask if the orbits of $W_p$ are
unions of blocks. We will show that this is the case, and give
examples in Section \ref{absect} to show that indeed these orbits are
not in general single blocks. (A similar result for the symplectic
Schur algebra has been given by the second author \cite{devblock}.)
Throughout the next two sections we will assume that we are working
over a field of characteristic $p>0$.

We will need the following positive characteristic analogue of
\cite[Proposition 4.2]{cdm}. Denote by $[\lambda]$ the set of boxes in
$\lambda$.

\begin{prop}\label{content}
 Let $\lambda, \mu\in X^+$ with $|\lambda|-|\mu|=2t\geq 0$.  If there
exists $M\leq \Delta_n(\mu^T)$ with
$$\Hom_{B_n(\delta)}(\Delta_n(\lambda^T),\Delta_n(\mu^T)/M)\neq 0$$ then
$$t(\delta-1)+\sum_{d\in[\lambda]}c(d)-\sum_{d\in[\mu]}c(d)\equiv 0 \mod p.$$
\end{prop}
\begin{proof}
As the localisation functor is exact, we may assume without loss of
generality that $\lambda\vdash n$. For $1\leq i<j\leq n$ let
$X_{i,j}$ be the Brauer diagram with edges between $l$ and $\bar{l}$
for all $l\neq i,j$, and with edges between $i$ and $j$ and $\bar{i}$
and $\bar{j}$. Then we define $T_n\in B_n(\delta)$ by
$$T_n=\sum_{1\leq i<j\leq n}X_{i,j}.$$ As in \cite[Lemma 4.1]{cdm} we
have for all $y\in\Delta_n(\mu^T)$ that
$$T_ny=
\left[t(\delta-1)-\sum_{d\in[\mu]}c(d)+
\sum_{1\leq i<j\leq n}(i,j)\right]y
$$ where $(i,j)$ denotes the transposition in $\Sigma_n$ permuting $i$
and $j$, regarded as an element of $B_n(\delta)$.

Let $\phi:\Delta_n(\lambda^T)\too\Delta_n(\mu^T)/M$ be a non-zero
$B_n(\delta)$-homomorphism. By our assumption on $n$, we have
$\Delta_n(\lambda^T)\cong S^{\lambda^T}$ (the Specht module labelled by
$\lambda^T$) as a module for $\Sigma_n$. As such modules are defined
over $\ZZ$, the remarks in the proof of \cite[Proposition 4.2]{cdm}
about the action of $\sum_{1\leq i<j\leq n}(i,j)$ on
$\Delta_n(\lambda^T)$ still hold, and hence this element acts as the
scalar $\sum_{d\in[\lambda]}c(d)$ on $\Delta_n(\lambda^T)$. Hence 
$$\sum_{1\leq i<j\leq
n}(i,j)\phi(x)=\sum_{d\in[\lambda]}c(d)\phi(x)$$ for all
$x\in\Delta_n(\lambda^T)$, and so for all $y+M\in\im(\phi)$ we must have
$$T_n(y+M)=
\left[t(\delta-1)-\sum_{d\in[\mu]}c(d)+
\sum_{d\in[\lambda]}c(d)\right](y+M).
$$ 

Again as in the proof of \cite[Proposition 4.2]{cdm}, the element
$T_n$ must act as zero on $\Delta_n(\lambda^T)$, and so 
$$t(\delta-1)-\sum_{d\in[\mu]}c(d)+
\sum_{d\in[\lambda]}c(d)\equiv 0\mod p$$
as required.
\end{proof}

We wish to replace the role played by the combinatorics of partitions
 by the action of our affine reflection group $W_p$. 

\begin{thm}
Suppose that $\lambda ,\mu\in X^+$. If there exists $M\leq
\Delta_n(\mu^T)$ with
$$\Hom_{B_n(\delta)}(\Delta_n(\lambda^T),\Delta_n(\mu^T)/M)\neq 0$$
then $\mu\in W_p.\lambda$.
\end{thm}

\begin{proof}
First note that $\Hom_{B_n(\delta)}(\Delta_n(\lambda^T),
\Delta_n(\mu^T)/M)\neq 0$ implies that $|\lambda | - |\mu| =2t\geq
0$. As if we had $|\lambda |<|\mu|$ then using the fact that the
localisation function $F$ is exact, we can assume that $\mu\vdash n$,
so $\Delta_n(\mu^T) \cong S^{\mu^T}$. However, this module only contains
composition factors of the form $L_n(\eta)$ where $\eta \vdash
n$, which gives a contradiction.

We now use induction on $n$. If $n=1$ then $\lambda = \mu = (1)$ and
so there is nothing to prove. Assume $n>1$. If $\lambda =\emptyset$
then by the above remark we have $\mu=\emptyset$ and we are done. Now
suppose that $\lambda$ has a removable box in row $i$ say. Then we
have
$${\rm Ind}\, \Delta_{n-1}((\lambda - \epsilon_i)^T)
\twoheadrightarrow \Delta_n(\lambda ^T)$$ and so, using our assumption
we have
$$\Hom_{B_n(\delta)}({\rm Ind}\, \Delta_{n-1}((\lambda -
\epsilon_i)^T), \Delta_n(\mu^T)/M)$$
$$=\Hom_{B_{n-1}(\delta)}(\Delta_{n-1}((\lambda - \epsilon_i)^T), {\rm
Res}\, (\Delta_n(\mu^T)/M)) \neq 0.$$
 Thus either (Case 1) we have
$$\Hom_{B_{n-1}(\delta)}(\Delta_{n-1}((\lambda - \epsilon_i)^T),
\Delta_{n-1}((\mu - \epsilon_j)^T)/N)\neq 0$$ for some positive
integer $j$ with $\mu - \epsilon_j\in X^+$ and some $N\leq
\Delta_{n-1}((\mu-\epsilon_j)^T)$,\\ 
or (Case 2) we have
$$\Hom_{B_{n-1}(\delta)}(\Delta_{n-1}((\lambda - \epsilon_i)^T),
\Delta_{n-1}((\mu + \epsilon_j)^T)/N)\neq 0$$ for some positive
integer $j$ with $\mu + \epsilon_j\in X^+$ and some $N\leq
\Delta_{n-1}((\mu+\epsilon_j)^T)$.

\noindent \textbf{Case 1:} Using Proposition \ref{content} for
$\lambda$ and $\mu$ and for $\lambda - \epsilon_i$ and $\mu -
\epsilon_j$, we see that $c(\lambda)_i \equiv c(\mu)_j \,\, \mbox{mod
$p$}$. Now, using induction on $n$ we have that $\mu - \epsilon_j \in
W_p \cdotp (\lambda - \epsilon_i)$. By Proposition \ref{orbits2}, we
can find $\pi\in \Sigma_n$ and $\sigma \, : \, {\bf n} \rightarrow
\{\pm 1\}$ such that $d(\sigma)$ is even and for all $1\leq m\leq n$,
if $\sigma (m)=1$ we have
$$c(\mu - \epsilon_i)_m \equiv c(\lambda)_{\pi(m)} \,\, \mbox{mod $p$}$$
and if $\sigma(m)=-1$ we have
$$c(\mu)_m + c(\lambda)_{\pi(m)} \equiv 2-\delta \,\, \mbox{mod
$p$}.$$ We will now construct $\pi '\in \Sigma_n$ and $\sigma ' \, :
\, {\bf n} \rightarrow \{ \pm 1\}$ to show that $\mu\in W_p \cdotp
\lambda$.  Suppose $\pi(j)=k$ and $\pi(l)=i$ for some $k,l\geq
1$. Define $\pi '$ by $\pi '(j)=i$, $\pi '(l)=k$ and $\pi '(m)=\pi(m)$
for all $m\neq j,l$. Now if $\sigma(j)=\sigma(l)$ then define $\sigma
'$ by $\sigma '(j)=\sigma '(l)=1$ and $\sigma '(m)=\sigma(m)$ for all
$m\neq j,l$. And if $\sigma(j)=-\sigma(l)$ the define $\sigma '$ by
$\sigma '(j)=1$, $\sigma '(l)=-1$ and $\sigma '(m)=\sigma (m)$ for all
$m\neq j,l$. Now it's easy to check, using the fact that
$c(\mu)_j\equiv c(\lambda)_i\,\, \mbox{mod $p$}$, that $\pi '$ and
$\sigma '$ satisfy the conditions in Proposition \ref{orbits2} for
$\lambda$ and $\mu$, and so
$\mu \in W_p \cdotp \lambda$.

\noindent \textbf{Case 2:} This case is similar to Case 1. Using
Proposition \ref{content} we see that $c(\lambda )_i + c(\mu)_j \equiv
2-\delta \,\, \mbox{mod $p$}$. Now using induction on $n$ we have
$\pi\in \Sigma_n$ and $\sigma \, : \, {\bf n} \rightarrow \{ \pm
1\}$ satisfying the conditions in Proposition \ref{orbits2} for
$\lambda -\epsilon_i$ and $\mu + \epsilon_j$.  Suppose $\pi(j)=k$ and
$\pi(l)=i$ for some $k,l\geq 1$. Define $\pi '$ by $\pi '(j)=i$, $\pi
'(l)=k$ and $\pi '(m)=\pi (m)$ for all $m\neq j,l$. Now if $\sigma
(j)=\sigma(l)$ the define $\sigma '$ by $\sigma '(j)=-1$, $\sigma
'(l)=-1$ and $\sigma '(m)=\sigma(m)$ for all $m\neq j,l$. And if
$\sigma (j)=-\sigma (l)$ then we define $\sigma '(j)=-1$, $\sigma
'(l)=1$ and $\sigma '(m)=\sigma (m)$ for all $m\neq j,l$.
\end{proof}

Note that by the cellularity of $B_n(\delta)$ this immediately implies

\begin{cor}
Two simple $B_n(\delta)$-modules $L(\lambda^T)$ and $L(\mu^T)$ are in
the same block only if $\mu\in W_p.\lambda$.
\end{cor}

Thus we have the desired necessary condition in terms of the affine
Weyl group for two weights to lie in the same block.

\section{Abacus notation and orbits of the affine Weyl group}
\label{absect}

In this section we will show that, even if $n$ is arbitrarily large,
being in the same orbit under the affine Weyl group is not
sufficient to ensure that two weights lie in the same block. This is
most conveniently demonstrated using the abacus notation \cite{jk},
and so we first explain how this can be applied in the Brauer algebra
case. We begin by recalling the standard procedure for constructing an
abacus from a partition, and then show how this is compatible with the
earlier orbit results for $W_p$. As in the preceding section, we
assume that our algebra is defined over some field of characteristic
$p>2$.

\begin{figure}[ht]
\includegraphics{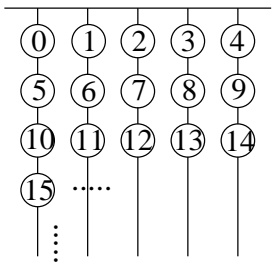}
\caption{The possible bead positions with $p=5$.}
\label{abacus}
\end{figure}

To each partition we shall associate a certain configuration of beads
on an abacus in the following manner. An {\it abacus with $p$ runners}
will consist of $p$ columns (called runners) together with some number
of beads distributed amongst these runners. Such beads will lie at a
fixed height on the abacus, and there may be spaces between beads on
the same runner. We will number the possible bead positions from left
to right in each row, starting from the top row and working down, as
illustrated in Figure \ref{abacus}.

\begin{figure}[ht]
\includegraphics{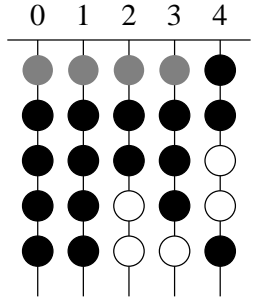}
\caption{The abacus for $(5,3,3,2,1,1,0^{10})$ when  
 $p=5$, $n=16$, and $b=20$.}
\label{exab}
\end{figure}

For a fixed value of $n$, we will associate to each partition
$\lambda=(\lambda_1,\lambda_2,\ldots,\lambda_t)$ of $m$, with $m\leq n$
and $n-m$ even, a configuration of beads on the abacus.  Let $b$ be a
positive integer such that $b\geq n$. We then represent $\lambda$ on
the abacus using $b$ beads by placing a bead in position numbered
$$\lambda_i+b-i$$ for each $1\leq i\leq b$, where we take
$\lambda_i=0$ for $i>t$. In representing such a configuration we will
denote the beads for $i\leq n$ by black circles, for $n<i\leq b$ by
grey beads, and the spaces by white circles (or blanks if this is
unambiguous). Runners will be numbered left to right from $0$ to
$p-1$. For example, the abacus corresponding to the partition
$(5,3,3,2,1,1,0^{10})$ when $p=5$, $n=16$, and $b=20$ is given in Figure
\ref{exab}. Note that the abacus uniquely determines the partition
$\lambda$.

We would like a way of identifying whether two partitions $\lambda$
and $\mu$ are in the same $W_p$ orbit directly from their abacus
representation. First let us rephrase the content condition which we
had earlier.

Recall from Proposition \ref{orbits2} and the definition of
$c(\lambda)$ that $\lambda$ and $\mu$ are in the same $W_p$-orbit if
and only if there exists $\pi\in\Sigma_n$ such that for each $1\leq
i\leq n$ either
$$\mu_i-i\equiv\lambda_{\pi(i)}-\pi(i)\mod p$$
or
$$\mu_i-i\equiv \delta-2-(\lambda_{\pi(i)}-\pi(i))\mod p$$ and
the second case occurs an even number of times.

Choose $b\in\NN$ with $2b\equiv 2-\delta\mod p$ (such a $b$ always
exists as $p>2$). Then $\lambda$ and $\mu$ are in the same $W_p$-orbit
if and only if there exists $\pi\in\Sigma_n$ such that for each $1\leq
i\leq n$ either
\begin{equation}\label{samerun}
\mu_i+b-i\equiv\lambda_{\pi(i)}+b-\pi(i)\mod p
\end{equation}
or
\begin{equation}\label{diffrun}
\mu_i+b -i\equiv p-(\lambda_{\pi(i)}+b -\pi(i))\mod p
\end{equation}
and the second case occurs an even number of times.  Thus if we also
choose $b$ large enough such that $\lambda$ and $\mu$ can be
represented on an abacus with $b$ beads then (\ref{samerun}) says that
the bead corresponding to $\mu_i$ is on the same runner as the bead
corresponding to $\lambda_{\pi(i)}$, and (\ref{diffrun}) says that the
bead corresponding to $\mu_i$ is on runner $l$ only if the bead
corresponding to $\lambda_{\pi(i)}$ is on runner $p-l$. Note that for
corresponding black beads on runner $0$ both (\ref{samerun}) and
(\ref{diffrun}) hold, and so we can use this pair of beads to modify
$d(\sigma)$ to ensure that it is even. Obviously if there are no such
black beads then the number of beads changing runners between
$\lambda$ and $\mu$ must be even.  Further, the grey beads (for
$i>n$) are the same on each abacus. Summarising, we have

\begin{prop}\label{abblock}
Choose $b\geq n$ with $2b\equiv 2-\delta\mod p$, and $\lambda$
and $\mu$ in $\Lambda_n$. Then
$\lambda$ and $\mu$ are in the same $W_p$-orbit if and only if

(i) the number of beads on runner $0$ is the same for $\lambda$ and
$\mu$, and

(ii) for each $1\leq l\leq p-1$, the \emph{total} number of beads
on runners $l$ and $p-l$ is the same for $\lambda$ and $\mu$, and

(iii) if there are no black beads on runner $0$, then 
the number of beads changing runners between $\lambda$ and $\mu$ must
be even.
\end{prop}

Note that condition (iii) plays no role when $n$ is large (compared to
$p$) as in such cases every partition will have a black bead on runner
$0$.

\begin{figure}[ht]
\includegraphics{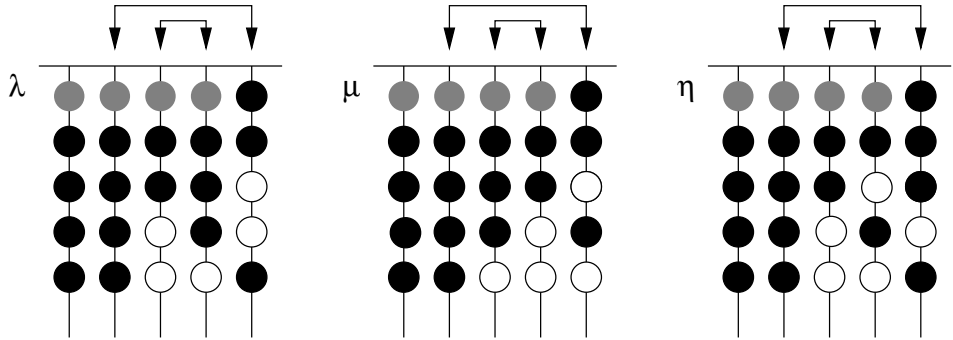}
\caption{Abacuses representing the elements $\lambda$, $\mu$ and
 $\eta$ in (\ref{trio}) with $b=20$.}
\label{threeab}
\end{figure}

To illustrate this result, consider the case $n=16$ and the partitions 
\begin{equation}\label{trio}
\lambda=(5,3,3,2,1,1),\quad \mu=(2,2,2,1,1,1),\quad
\eta=(5,3,3,2,1,1,1).
\end{equation}
Take $p=5$ and $\delta=2$, then $b=20$
satisfies $2b\equiv 2-\delta\mod p$, and is large enough for all three
partitions to be represented using $b$ beads. The respective abacuses
are illustrated in Figure \ref{threeab}, with the matching rows for
condition (ii) in Proposition \ref{abblock} indicated.

We see that $\mu\in W_p.\lambda$, as the number of beads on runner
$0$, and on runners $1/4$ and $2/3$ are the same for both $\lambda$
and $\mu$ (respectively 5, 8, and 7) and there is a black bead on
runner $0$. (The number of beads moving from runner $l$ to a distinct
runner $p-l$ is $1$, which is odd. However, as discussed above, we can
chose $\sigma$ such that one of the two black beads on runner $0$ is
regarded as moving (to the same runner), to obtain the required even
number of such moves. If there were no black beads on runner $0$ then
this would not be possible.) However, $\eta\notin W_p.\lambda$ as
columns $1/4$ and $2/3$ have 9 and 6 entries respectively.

Having reinterpreted the orbit condition in terms of the abacus, we
will now show that the orbits of $W_p$ can be
\emph{non-trivial} unions of blocks for $B_n(\delta)$.

\begin{thm} 
Suppose that $k$ is of characteristic $p>2$. Then for arbitrarily large
$n$ there exist $\lambda\vdash n$ and $\mu\vdash n-2$ which are in the
same $W_p$-orbit but not in the same $B_n(\delta)$-block.
\end{thm}
\begin{proof}
Let $b\in\NN$ be such that
$$2b\equiv 2-\delta\mod p.$$ If $b$ is even (respectively odd),
consider the partial abacuses illustrated in Figure \ref{even}
(respectively Figure \ref{odd}).

\begin{figure}[ht]
\includegraphics{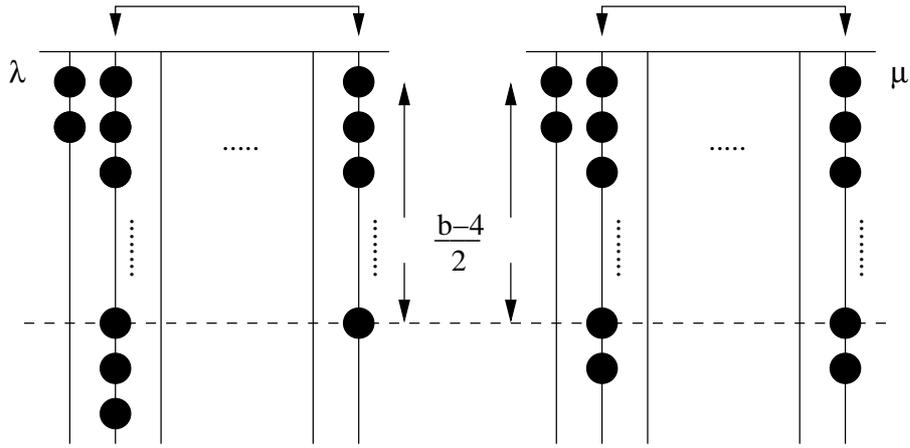}
\caption{The partial abacuses for $\lambda$ and $\mu$ when $b$ is even.}
\label{even}
\end{figure} 

\begin{figure}[ht]
\includegraphics{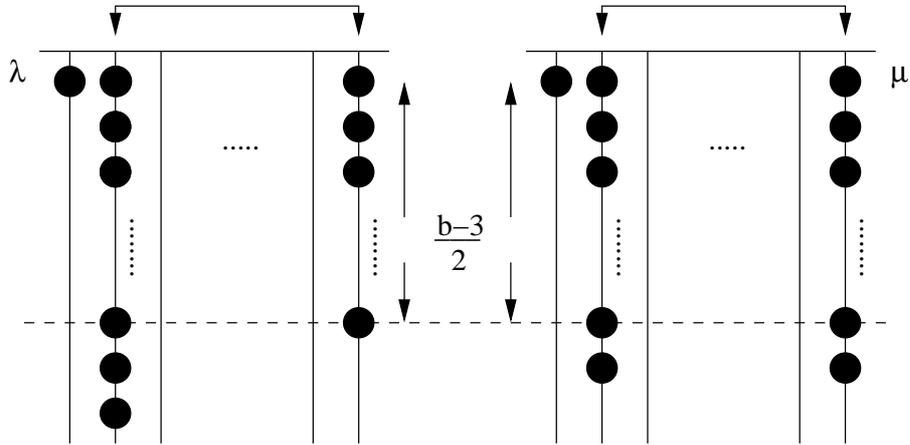}
\caption{The partial abacuses for $\lambda$ and $\mu$ when $b$ is
odd.}
\label{odd}
\end{figure} 

These will not correspond directly to partitions $\lambda$ and $\mu$,
as the degree of each partition will be much larger than $b$. However,
by completing each in the same way (by adding the same number of black
beads in rows from right to left above each partition, followed by a
suitable number of grey beads), they can be adapted to form abacuses
of partitions $\lambda\vdash n$ and $\mu\vdash n-2$ for some $n>>0$
and for some $b'\equiv b\mod p$. (This corresponds to adding
sufficiently many zeros to the end of each partition such that each
has $|\lambda|$ parts.)

It is clear from Proposition \ref{abblock} that in each case $\lambda$
and $\mu$ are in the same $W_p$-orbit. Note that for both $\lambda$
and $\mu$, all beads are as high as they can be on their given runner.
If we move any bead to a higher numbered position then this
corresponds to increasing the degree of the associated partition.
Thus $\lambda$ and $\mu$ are the only partitions with degree at most
$|\lambda|$ in their $W_p$-orbit.Also it is easy to check that $\mu$
is obtained from $\lambda$ by removing two boxes from the same
row. Clearly by increasing $b$ we can make $n$ arbitrarily large.

To complete the proof, it is enough to show that $L_n(\lambda^T)$ and
$L_n(\mu^T)$ are not in the same $B_n(\delta)$-block. We will reduce
this to a calculation for the symmetric group, and use the
corresponding (known) block result in that case. To state this we need
to recall the notion of a $p$-core. 

A partition is a \emph{$p$-core} if the associated abacus has no gap
between any pair of beads on the same runner. We associate a unique
$p$-core to a given partition $\tau$ by sliding all beads in some
abacus representation of $\tau$ as far up each runner as they can go,
and taking the corresponding partition.  By the Nakayama conjecture
(see \cite{meiertappe} for a survey of its various proofs), two
partitions $\tau$ and $\eta$ are in the same block for $k\Sigma_n$ if
and only if they have the same $p$-core. It is also easy to show
(using the definition of $p$-cores involving the removal of $p$-hooks
\cite{mathas}) that if $\tau$ is a $p$-core then so is $\tau^T$.

Returning to our proof, as $\lambda\vdash n$ we have that the cell
module $\Delta_n(\lambda^T)$ is isomorphic to the Specht module
$S^{\lambda^T}$ as a $k\Sigma_n$-module (by \cite[Section 2]{dhw}). As
$\lambda$ is a $p$-core so is $\lambda^T$, and hence
$\Delta_n(\lambda^T)$ is in a $k\Sigma_n$-block on its own (by the
Nakayama conjecture) so is simple as a $k\Sigma_n$-module,
isomorphic to $D^{\lambda^T}$, and hence equal to $L_n(\lambda^T)$ as a
$B_n(\delta)$-module.

If 
$$[\Delta_n(\mu^T):L_n(\lambda^T)]\neq 0$$
then we must have
$$[\res_{k\Sigma_n}\Delta_n(\mu^T):D^{\lambda^T}]\neq 0.$$ However,
$\res_{k\Sigma_n}\Delta_n(\mu^T)$ has a Specht filtration where the
multiplicity of $S^{\eta^T}$ in this filtration is given by the
Littlewood-Richardson coefficient $c_{\mu^T(2)}^{\eta^T}$
\cite[Proposition 8]{pagetLRR}. 
In particular, as $\lambda^T$ is obtained from $\mu^T$ by adding two
boxes in the same column we see that $S^{\lambda^T}=D^{\lambda^T}$
does not appear as a Specht subquotient in this filtration
\cite[remarks after Theorem 3.1]{dhw}. However, we still have to prove
that it cannot appear as a composition factor of some other
$S^{\eta^T}$. But this is clear, as if it did then $\eta^T$ would have
to have the same $p$-core as $\lambda$, but $\lambda$ is already a
$p$-core and hence this is impossible.

This proves that $\Delta_n(\mu^T)=L_n(\mu^T)$, and so $\lambda$ and $\mu$
are in different blocks for $B_n(\delta)$.
\end{proof}

\begin{figure}[ht]
\includegraphics{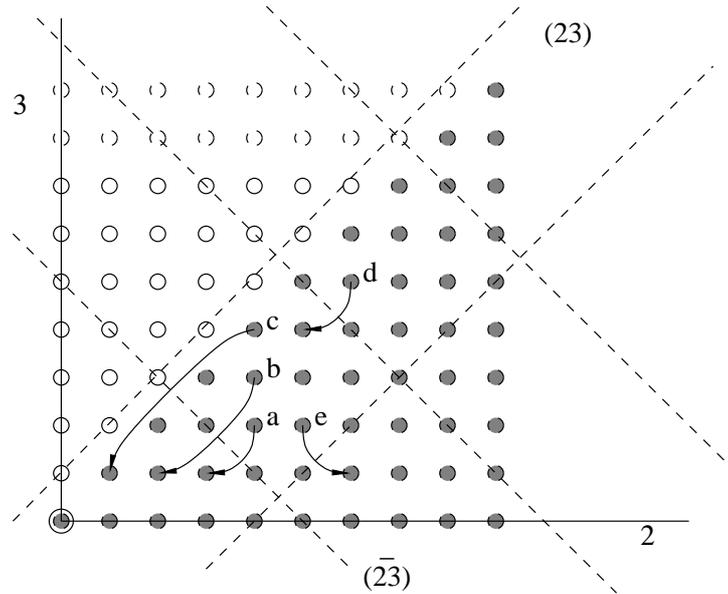}
\caption{Assorted examples 
  with $p=5$ and $\delta=2$.}
\label{rosetta}
\end{figure}

\begin{figure}[ht]
\includegraphics{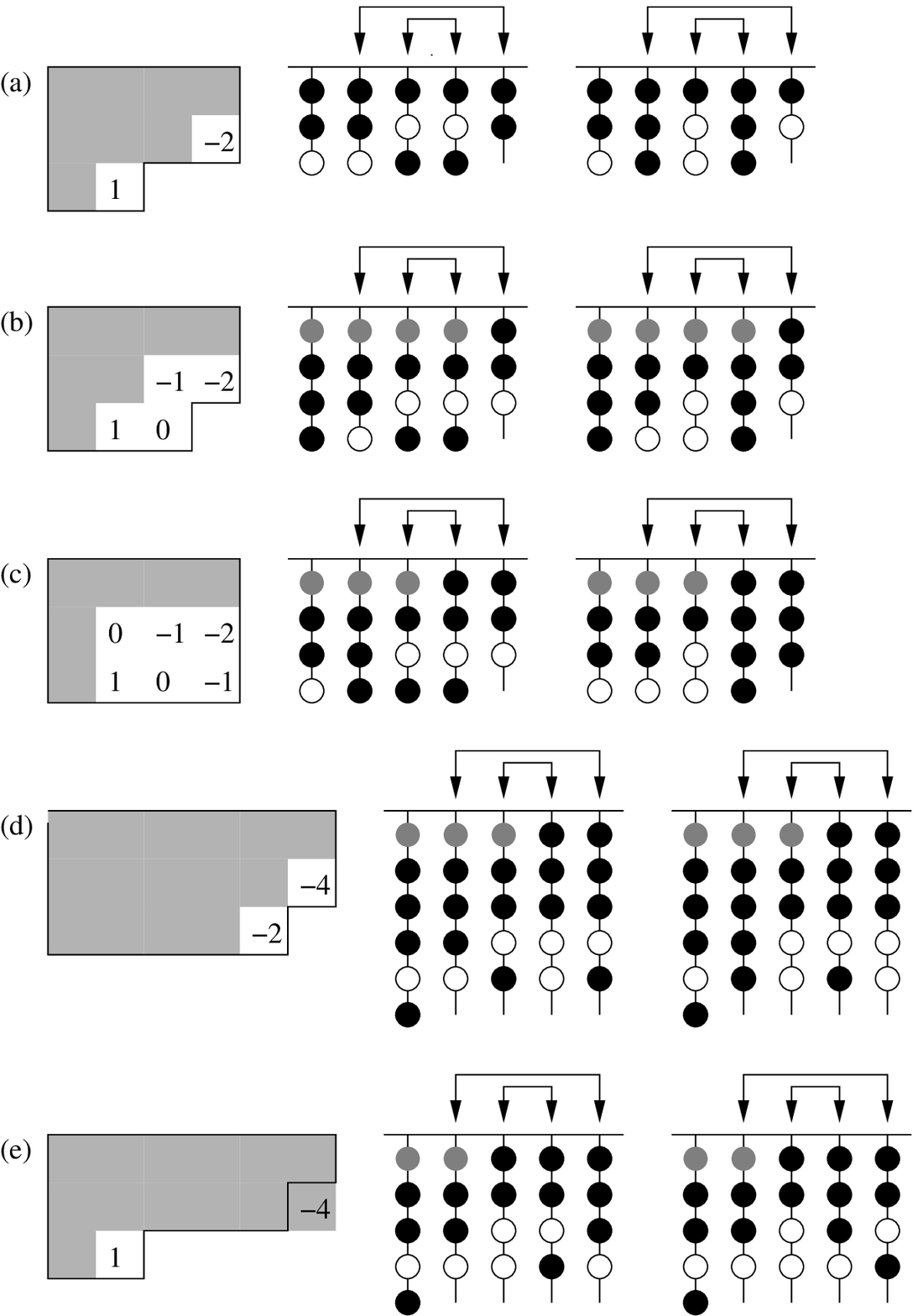}
\caption{Assorted examples 
  with $p=5$ and $\delta=2$.}
\label{rosetta2}
\end{figure}

To conclude, we illustrate some examples of various affine reflections
together with the corresponding partitions and abacuses, when $p=5$
and $\delta=2$. Our condition on $b$ implies that it must be chosen to
be a multiple of $5$. Reflections are labelled (a)--(e) in Figure
\ref{rosetta}, with the corresponding partitions and abacuses in
Figure \ref{rosetta2}.  Case (a) corresponds to the reflection from
$(4,4,2)$ to $(4,3,1)$, with $n=b=10$. Case (b) corresponds to the
reflection from $(4,4,3)$ to $(4,2,1)$, with $n=11$ and $b=15$.  Case
(c) corresponds to the reflection from $(4,4,4)$ to $(4,1,1)$ with
$n=12$ and $b=15$. These three cases only use elements from $W$, and
so would be reflections in any characteristic. Hence the condition on
matched contents in these cases are equalities, not merely
equivalences mod $p$. Case (d) corresponds to the reflection from
$(6,6,5)$ to $(6,5,4)$ with $n=17$ and $b=20$. This is a strictly
affine phenomenon, and so the paired boxes only sum to $1-\delta$ mod
$p$. Finally, case (e) corresponds to the reflection from $(6,5,2)$ to
$(6,6,1)$ with $n=13$ and $b=15$. This is our only example of
reflection about an affine $(ij)$ line, and so is the only case
illustrated where the number of boxes is left unchanged.


\providecommand{\bysame}{\leavevmode\hbox to3em{\hrulefill}\thinspace}
\providecommand{\MR}{\relax\ifhmode\unskip\space\fi MR }
\providecommand{\MRhref}[2]{%
  \href{http://www.ams.org/mathscinet-getitem?mr=#1}{#2}
}
\providecommand{\href}[2]{#2}

\end{document}